\documentclass[12pt]{article}
\usepackage{amsxtra,amssymb,amsthm,amsmath,latexsym}

\theoremstyle{plain}
\newtheorem{theorem}{Theorem}

\newtheorem{remark}{Remark}

\numberwithin{equation}{section}

\newcommand{\refT}[1]{Theorem~\ref{T:#1}}
\newcommand{\refS}[1]{Section~\ref{S:#1}}

\def\ds{\displaystyle}

\def\bysame{\rule{.5in}{.005in},\ }

\def\calE{{\mathcal E}}

\def\dotu{\dot{u}}
\def\dotW{\dot{W}}

\def\oH{{\overset{\circ}{H}}}
\def\oH1{{\overset{\circ}{H}\kern-.02in{}^1}}

\def\const{{\hbox{\,const\,}}}

\def\bee{\begin{equation*}}
\def\eee{\end{equation*}}
\def\be{\begin{equation}}
\def\ee{\end{equation}}

\begin{document} 
\title{ Dynamical systems
method for solving linear ill-posed problems }

\author{A.G. Ramm\\
 Mathematics Department, Kansas State University, \\
 Manhattan, KS 66506-2602, USA\\
ramm@math.ksu.edu,\\ fax 785-532-0546, tel. 785-532-0580}

\date{}
\maketitle\thispagestyle{empty}

\begin{abstract}
\footnote{key words: ill-posed problems, dynamical systems method (DSM),
regularization parameter, discrepancy principle, unbounded operators,
linear operator equations   
}
Various versions of the Dynamical Systems Method (DSM) 
are proposed for solving linear ill-posed problems
with bounded and unbounded operators. Convergence of the proposed methods
is proved. 
Some new results concerning discrepancy principle for choosing 
regularization parameter are obtained.

\end{abstract}

\section{Introduction}\label{S:1}
In this paper part of the results from the author's invited 
talk at ICMAAS06 are presented.
This part deals with linear ill-posed problems.

Some of the results in this paper are taken from the 
papers of the 
author cited in the bibliography and many results are new. 
The DSM
method is developed in \cite{R}, \cite{R1}-\cite{R21}. The discrepancy 
principle, discussed earlier in 
\cite{M},  has been studied recently in \cite{R1}-\cite{R3} and in 
\cite{R}.

Consider an equation
\be\label{e1.1} Au-f=0, \ee
where $A$ is an operator in a Banach space $X$. If $A$ is a homeomorphism of $X$ onto $X$ (i.e., 
continuous injective and surjective map in $X$ which has a continuous inverse map) 
then problem \eqref{e1.1} is called well-posed in the Hadamard's sense.
Otherwise it is called ill-posed.

The Dynamical Systems Method (DSM) for solving equation \eqref{e1.1} consists of solving the Cauchy
problem
\be\label{e1.2}
\dotu=\Phi(t,u),\qquad u(0)=u_0;\qquad \dotu=\frac{du}{dt}, 
\ee
which we call a dynamical system, where $\Phi$  is chosen so 
that
\be\label{e1.3}
\exists ! u(t)\quad \forall t\geq 0;\qquad \exists u(\infty);\qquad A(u(\infty))=f.\ee
Here $u(\infty):=\lim_{t\to\infty}u(t)$.
We do not assume that the solution to \eqref{e1.1} is unique, but we do assume that there exists
a solution to \eqref{e1.1}. If \eqref{e1.3} holds, then we say that DSM is justified for solving
equation \eqref{e1.1}.

There is a large body of literature on solving ill-posed problems 
(see e.g., \cite{I}, \cite{R}, \cite{T}
and references therein). Variational regularization, iterative 
regularization, quasisolutions and quasiinversion 
are some of the methods for stable solution of ill-posed 
problems discussed in the literature.
In this paper several new methods for stable solutions of 
linear ill-posed problems is discussed. They are based on the dynamical 
systems method (DSM).
We describe the DSM for solving ill-posed problems.
This method was developed relatively recently as a method for 
solving a wide variety of linear and 
nonlinear ill-posed problems \cite{R1}-\cite{R21} although it 
was proposed already in \cite{G} 
for solving well-posed problems.

Because of the space limitations we will not discuss solving nonlinear ill-posed problems 
by the DSM in this paper, and refer the reader to \cite{R}, \cite{R9}.

Here we describe a new version of DSM for solving linear ill-posed problems. There are many practical
problems of this type. We only mentioned solving ill-conditioned linear algebraic systems and 
Fredholm equations of the first kind. The novel points in our results are not only the method of solving these 
problems by DSM but also the applicability of the method to unbounded operators. In the literature
a widely-discussed method for solving ill-posed problem \eqref{e1.1} is the method of variational 
regularization  introduced by Phillips (\cite{P}) and studied by 
Tikhonov (\cite{T}), Morozov (\cite{M}), Ivanov (\cite{I}), and many other 
authors
under the assumption that the operator $A$ in \eqref{e1.1} 
is a linear bounded operator. There are some results on
regularization of unbounded operators also (e.g., \cite{L}, 
\cite{R21}).
The variational regularization method for stable solution of \eqref{e1.1} consists of
solving the problem
\be\label{e1.4} F(u):=\|Au-f_\delta\|^2+a\|u\|^2=\min, 
\ee
where $a>0$ is a constant, called a regularization parameter, $f_\delta$ are the ``noisy data", 
i.e.,  an element which satisfies the inequality $\|f_\delta-f\|\leq\delta$ and which is given
together with the ``noise level" $\delta>0$, while the exact data $f$ are not known. The stable solution to 
\eqref{e1.1} is an element $u_\delta$ such that $\lim_{\delta\to 0}\|u_\delta-y\|=0$,
where $Ay=f$ and $y$ is the unique minimal-norm solution of the linear equation \eqref{e1.1}.
If $X$ is a Hilbert space $H$, which we assume below, then the minimal-norm solution is the 
solution which is orthogonal to the null space $N$ of $A$, $N=N(A)=\{u:Au=0\}$.
If the linear operator $A$ in \eqref{e1.1} is unbounded, then we assume that it is closed and
densely defined, so that its adjoint $A^\ast$ is densely defined and closed
(see, e.g.,  \cite{K}).

If $A$ is bounded, then the necessary and sufficient condition for the minimizer of the quadratic 
functional \eqref{e1.4} is the Euler's equation 
\be\label{e1.5} T_a u=A^\ast f_\delta,\quad T_a:=T+aI,\quad 
T:=A^\ast A, \ee
where $I$ is the identity operator and $T\geq 0$ is a 
selfadjoint operator. Equation \eqref{e1.5} 
has a unique solution  $u_{a,\delta}=T^{-1}_a A^\ast f_\delta$.
One can choose $a=a(\delta)$ so that $\lim_{\delta\to 0}a(\delta)=0$
and $u_\delta:=u_{a(\delta),\delta}$ is a stable solution to 
\eqref{e1.1}:
\be\label{e1.6} \lim_{\delta\to 0} \|u_\delta-y\|=0, \ee
where $Ay=0$, $y\perp N$. There are a priori choices for choosing $a(\delta)$
and a posteriori ones. An a priori choice is based on the estimate:
\be\label{e1.7}
 \|T^{-1}_a A^\ast f_\delta-y\| \leq \|T^{-1}_a 
A^\ast(f_\delta-f)\|
+\|T^{-1}_a A^\ast f-y\| \leq \frac{\delta}{2\sqrt{a}} 
+\eta(a), \ee
where
\be\label{e1.8}
\eta^2(a)=\|T^{-1}_a Ty-y\|^2 = a^2 \int_0^\infty  
\frac{d(E_s y,y)}{(s+a)^2}
\to \|P_Ny\|^2 \hbox{\ as\ } a\to 0,\ee
and $P_N$ is the orthoprojector onto $N$. Since we assume that $y\perp N$,
equation \eqref{e1.8} implies
\be\label{e1.9} \lim_{a\to 0} \eta(a)=0.\ee
The term $\frac{\delta}{2\sqrt{a}}$ in \eqref{e1.7} appears due to the 
estimate:
\be\label{e1.10}
\|T^{-1}_a A^\ast\| = \|A^\ast Q^{-1}_a\| 
= \|U Q^{\frac{1}{2}} Q^{-1}_a\|\leq \|Q^{\frac{1}{2}} Q^{-1}_a\|
=\sup_{s\geq 0}\frac{\sqrt{s}}{s+a} = \frac{1}{2\sqrt{a}}. \ee
Here $Q:=AA^\ast$, $U$ is a partial isometry, $Q_a:=Q+aI$,
we have used the formula
\be\label{e1.11} T^{-1}_a A^\ast=A^\ast Q^{-1}_a,\ee
the polar decomposition: $A^\ast=U Q^{\frac{1}{2}}$, and the 
spectral theorem for the selfadjoint
operator $Q$, namely $\|g(Q)\|=\sup_{s\geq 0} |g(s)|$.
Formula \eqref{e1.11} is obvious if $A$ is bounded: multiply \eqref{e1.11} by $T_a$ from the left
and then by $Q_a$ from the right, and get $A^\ast(AA^\ast+aI)=(A^\ast A+aI)A^\ast$, 
which is an obvious identity. Since the operators $Q_a$ and $T_a$ are boundedly invertible, 
one may reverse the steps and get  \eqref{e1.11}. Thus a priori choices of 
$a(\delta)$, which imply \eqref{e1.6}, are:
\be\label{e1.12}
\lim_{\delta\to 0} \frac{\delta}{\sqrt{a(\delta)}}=0, \qquad \lim_{\delta\to 0} a(\delta)=0.\ee
There are many functions $a(\delta)$ satisfying \eqref{e1.12}.
One can find an optimal value $a(\delta)$ by minimizing the right-hand 
side of \eqref{e1.7} with respect to $a$. Alternatively, one may calculate
$a(\delta)$ by solving the equation $\delta= 2 \sqrt{a(\delta)}\eta(a)$
for $a$ for a fixed small $\delta >0$.

If $A$ is closed, densely defined in $H$, unbounded, and $a=\const>0$, 
then the author has proved that the operator 
$T^{-1}_aA^\ast$, with the domain
$D(A^\ast)$, is closable, its closure, denoted again $T^{-1}_a A^\ast$,
is a bounded operator defined on all of $H$, $\|T^{-1}_a A^\ast\|\leq\frac{1}{2\sqrt{a}}$,
and \eqref{e1.1} holds.
For convenience of the reader let us sketch the proof of these claims. To check that 
$T^{-1}_a A^\ast$ is closable, one takes $h_n\in D(T^{-1}_a A^\ast)=D(A^\ast)$,
such that $h_n\to 0$ and $T^{-1}_a A^\ast h_n\to g$ as $n\to\infty$, and
checks that $g=0$. Indeed, let $u\in H$ be arbitrary. Then
\be\label{e1.13}
(g,u)=\lim_{n\to\infty} (T^{-1}_a A^\ast h_n,u) = \lim_{n\to\infty} (h_n,A T^{-1}_au)=0. \ee
Since $u$ is arbitrary, this implies $g=0$, as claimed.
Note that $T^{-1}_a u\in D(A)$, so that the above calculation is justified. 
If one drops the index $n$ and the $\lim$ in \eqref{e1.13}, then one can see that the adjoint
to the closure of $T^{-1}_a A^\ast$ is the operator $A T^{-1}_a$, defined on all of $H$
and bounded
\bee \|A T^{-1}_a\| = \|U T^{\frac{1}{2}}T^{-1}_a\| \leq\frac{1}{2\sqrt{a}}. \eee
Since $\|A^\ast\|=\|A\|$, one gets 
$\|T^{-1}_a A^\ast\|\leq \frac{1}{2\sqrt{a}}$.
Finally formula \eqref{e1.11} can be proved for unbounded closed, densely define operator
$A$ as above, if one checks that the operator $A^\ast AA^\ast$  is densely defined.
This operator is densely defined because the operator $A^\ast AA^\ast A=T^2$ is densely
defined if $T$ is, and $D(T^2)\subset D(A^\ast AA^\ast)$.

Let us now describe an a posteriori choice of $a(\delta)$ 
which implies \eqref{e1.6} 
and which is called a discrepancy principle. 
This principle was discussed in 
\cite{M}, \cite{R}. It consists in finding 
$a(\delta)$ from the equation
\be\label{e1.14} \|A u_{a,\delta} -f_\delta\| = 
C\delta, \qquad 1<C<2, \ee
where $C=const$,  $\|F_\delta\|>C\delta$ and $u_{a,\delta} 
=T^{-1}_a A^\ast f_\delta$.
One can prove (e.e. see \cite{R}) that the equation \eqref{e1.14} for a small fixed 
$\delta>0$ has a unique solution $a(\delta)$, $\lim_{\delta\to 0}a(\delta)=0$, and 
$u_\delta=u_{a(\delta),\delta}$ satisfies \eqref{e1.6}, i.e. $u_\delta$ is a stable solution
to \eqref{e1.1}. 
To prove these claims one writes \eqref{e1.14} as
\be\label{e1.15}
C^2\delta^2= \|[A T^{-1}_a A^\ast-I ]f_\delta \|^2
=\|[QQ^{-1}_a-I]f_\delta\|^2
=\int^\infty_0 \frac{a^2d(E_\delta 
f_\delta,f_\delta)}{(s+a)^2} : = h(a,\delta), \ee
and takes into account that $h(a,\delta)$  is, for a fixed $\delta>0$,
a continuous monotone function of $a$, $h(\infty,\delta)=\|f_\delta\|^2>
C^2\delta^2$
and $h(+0,\delta)=\|P_{N^\ast}f_\delta\|\leq\delta^2$,
so that there exists a unique $a=a(\delta)$ such that $h(a(\delta),\delta)=C^2\delta^2$.
Here we have denoted $N^\ast:=N(A^\ast)$, used the obvious relation 
$N(Q)=N(A^\ast)$,
the inequality 
$\|P_{N^\ast}f_\delta\|\leq \|P_{N^\ast}(f_\delta-f)\| 
+ \|P_{N^\ast}f\| \leq \|f_\delta-f\|=\delta,$ and the relation $P_{N^\ast}f=0$.
This last relation follows from the formulas $f\in R(A)$ and $R(A)\perp N^\ast$.

Let us now check that if $a(\delta)$ solves \eqref{e1.5} then $u_\delta=a_{a(\delta),\delta}$ 
satisfies \eqref{e1.6}. One has $F(u_\delta)\leq F(y)$, so
\be\label{e1.16}
\|Au_\delta-f_\delta\|^2+a(\delta) \|u_\delta\|^2 \leq \delta^2+a(\delta) \|y\|^2. \ee
Since $\|Au_\delta-f_\delta\|^2 = C^2\delta^2>\delta^2$ and $a(\delta)>0$,
one gets
\be\label{e1.17} \|u_\delta\| \leq \|y\|. \ee
Therefore one can select a weakly convergent sequence $u_n=u_{\delta_n}$,
$u_n\rightharpoonup u$, as $n\to \infty$. 
Let us prove that $u=y$ and $\lim_{n\to\infty} \|u_n-y\|=0$.
Since this holds for any subsequence, it follows that \eqref{e1.6} holds.
To prove that $u=y$ note that \eqref{e1.17} implies $\|u\|\leq \|y\|$,
and that $u$ solves \eqref{e1.1}. Since the minimal-norm solution to \eqref{e1.1}
is unique, it follows that $u=y$. To check that $u$  solves \eqref{e1.1} we 
note that
$\lim_{\delta\to 0} \|Au_\delta-f\|=0$ as follows from \eqref{e1.16} because
$\lim_{\delta\to 0}a(\delta)=0$. The relations $u_\delta\rightharpoonup u$ 
and $\|Au_\delta-f\|\to 0$  as $\delta\to 0$
imply $Au=f$ and $\lim_{\delta\to 0}\|u_\delta-u\|=0$.
Let us first check that $Au=f$. One has 
\bee (f,g)=\lim_{\delta\to 0}(Au_\delta,g)=(u,A^\ast g)\qquad \forall g\in D(A^\ast).\eee
Thus $u\in D(A)$ and $Au=f$, as claimed.
As we proved above, this implies that $u=y$. Therefore $u_\delta\rightharpoonup y$,
$\|u_\delta\|\leq\|y\|$. This implies that 
$\lim_{\delta\to 0} \|u_\delta-y\|=0$. 
Indeed,
\be\label{e1.18}
\|u_\delta-y\|^2 = \|u_\delta\|^2 +\|u\|^2 -2 Re (u_\delta,y)
\leq 2\|y\|^2 - 2 Re(u_\delta,y) \to 0 \quad\hbox{\ as\ } \delta\to 0. \ee
Thus the relation \eqref{e1.6} is proved for the choice of $a(\delta)$ by the discrepancy principle.

The drawback of the a priori choice if $a(\delta)$ is that it is nonunique and although it
guarantees convergence \eqref{e1.16},  the error of the method can be large if 
$\delta>0$ is fixed. The drawback of the discrepancy principle is the necessity of solving
the nonlinear equation \eqref{e1.14} and also a possible large error for a fixed $\delta$.

In Section 2 we discuss the DSM for solving linear equations \eqref{e1.1}.

\section{DSM for solving linear problems}\label{S:2}

We assume first that the linear closed densely defined in $H$ operator in \eqref{e1.1} is selfadjoint,
$A=A^\ast$. This is not a restriction: every solvable linear equation \eqref{e1.1} is equivalent 
to an equation $Tu=A^\ast f$, where $T=T^\ast=A^\ast A$. Indeed, if $Au=f$, 
then applying $A^\ast$
and assuming $f\in D(A^\ast)$, one gets $Tu=A^\ast f$. 
Conversely, if $Tu=A^\ast f$, and $f=Ay$, then $Tu=Ty$. 
Multiply the equation $0=T(u-y)$ by $u-y$ and get 
$0=(A^\ast A(u-y),u-y)=\|Au-Ay\|^2$. Thus $Au=Ay=f$.
If $A$ is bounded, then $f\in D(A^\ast)$ for any $f\in H$.
If $A$ is unbounded, then $ D(A^\ast)$ is a linear 
dense subset of $H$. In this case, 
if $f\not\in D(A^\ast)$, then we define the solution of equation 
$Tu=A^\ast f$ by the formula $u=\lim_{a\to 0}T^{-1}_a A^\ast f$.
As we have proved in \refS{1}, for any $f\in R(A)$ this limit exists and equals
to the minimal-norm solution $y$: $\lim_{a\to 0} T^{-1}_a A^\ast Au=y$ 
if $Au=f$. This is true because $\lim_{a\to 0} T^{-1}_a Tu=u-P_Nu=y$.

The DSM for solving equation \eqref{e1.1} with a linear selfadjoint operator 
can be constructed as follows. Consider the problem 
\be\label{e2.1}
\dotu_a=i(A+ia)u_a-if, \qquad u(0)=0; \qquad \dotu=\frac{du}{dt}, \ee
where $a=const>0$. Our first result is formulated as \refT{1}.

\begin{theorem}\label{T:1}
If $Ay=f$ and $y\perp N$, then
\be\label{e2.2} \lim_{a\to 0}\  \lim_{t\to\infty} \  u_a(t)=y.\ee
\end{theorem}

Our second result shows that the method, based on Theorem 1, gives a stable 
solution of the equation $Au=f$.
Assume that $\|f_\delta-f\|\leq\delta$, and let $u_{a,\delta}(t)$ be the
solution to \eqref{e2.1} with $f_\delta$ in place of $f$.

\begin{theorem}\label{T:2}
There exist $t=t_\delta$, $\lim_{\delta\to 0} t_\delta=\infty$, and 
$a=a(\delta)$,  $\lim_{\delta\to 0} a(\delta)=0$,  such that 
$u_\delta:=u_{a(\delta),\delta}(t_\delta)$  satisfies \eqref{e1.6}.
\end{theorem}

We will discuss the ways to choose $a(\delta)$ and $t_\delta$ after the proofs
of these theorems are given.

From the numerical point of view if one integrates problem \eqref{e2.1} with the
exact data $f$ on the interval $0\leq t\leq T$, and $T$ is fixed, then one is
interested in choosing such $a=a(T)$ that 
$\lim_{T\to \infty} \|u_{a(T)} (T)-y\|=0$. We will give such a choice of $a(T)$.

Before we start proving these two theorems, let us explain the ideas of the proof. 
Suppose $B$ is a linear operator and its inverse  $B^{-1}$ exists and is bounded.

Finally, assume that
\be\label{e2.3} \lim_{t\to\infty} \|e^{Bt}\|=0. \ee
This will happen, for example, if $Re\,B\leq-c$, $c>0$ is a constant. Under these 
assumptions  one has
\be\label{e2.4}
\int^t_0 e^{Bs}ds=B^{-1} \left(e^{Bt}-I\right),
\quad -\lim_{t\to\infty} \int^t_0 e^{Bs}ds=B^{-1}. \ee
The operator $\int^t_0 e^{Bs}ds$ solves the problem
\be\label{e2.5}
  \dotW=BW+I, \qquad W(0)=0, \ee
where $I$ is the identity operator. If $\lim_{t\to\infty} \|e^{Bt}\|=0$,
then 
\be\label{e2.6}
  -\lim_{t\to\infty} \quad W(t)=B^{-1}. \ee

The basic idea of the DSM is the representation of the inverse operator 
as the limit as  $t\to\infty$ of the solution to the  Cauchy problem \eqref{e2.5}.

Let us prove \refT{1}.

\begin{proof}[Proof of \refT{1}.]
The solution to \eqref{e2.1} is 
\bee
u_a(t)=\int^t_0 e^{i(A+ia)(t-s)} (if)ds=[i(A+ia)]^{-1}
\left( e^{i(A+ia)t}-I\right)(-if).\eee
Since $\|e^{i(A+ia)t}\|=e^{-at}\to 0$ as $t\to\infty$, one 
gets
\be\label{e2.7} 
 \lim_{t\to\infty} u_a(t)=(A+ia)^{-1}f. \ee
Since $f=Ay$, one has
\be\label{e2.8}
 \eta(a):=\|(A+ia)^{-1} Ay-y\| =a \|(A+ia)^{-1}y\|\to 0
 \quad \hbox{\ as\ }a\to 0.\ee

 The last relation follows from the assumption $y\perp N$. 
Indeed, by the spectral theorem one has:
 \bee
 \begin{aligned}
 \lim_{a\to 0} \eta^2(a) 
 &= \lim_{a\to 0} a^2 \|(A+ia)^{-1}y\|^2\\
 &=\lim_{a\to 0} \int^\infty_{-\infty}\ \frac{a^2}{s^2+a^2}d (E_sy,y)
 =\|(E_0-E_{0-0})y\|^2
 =\|P_Ny\|^2=0. 
 \end{aligned}\eee
\refT{1} is proved.
\end{proof}

\begin{remark}\label{R:1}
In a numerical implementation of \refT{1} one chooses $\tau$, and $a=a(\tau)$, 
and integrates \eqref{e2.1} on the interval $[0,\tau]$.
One chooses $\tau$ so that
\be\label{e2.9}  
\lim_{\tau\to\infty} \|u_{a(\tau)}-y\|=0.\ee
In order to choose $a(\tau)$, note that  $\|u_a(t)-u_a(\infty)\|\leq \frac{e^{-at}}{a}$,
as follows from the derivation of \eqref{e2.7}. Therefore, by \refT{1}, the relation \eqref{e2.9}
holds if
\be\label{e2.10}
\lim_{\tau\to\infty} \frac{e^{-a(\tau)\tau}}{a(\tau)}=0\qquad
\lim_{\tau\to\infty} a(\tau)=0.\ee
For example, one may take $a(\tau)=\tau^{-\gamma}$, where $0<\gamma<1$ is a constant.
\end{remark}

\begin{proof}[Proof of Theorem 2.] 
Let us start with the formula:
\be\label{e2.11}
u_{a,\delta}(t)=\int^t_0 e^{i(A+ia)(t-s)} (-if_\delta)ds
=[i(A+ia)]^{-1} \left(e^{i(A+ia)t} -I\right) (if_\delta). \ee
Thus
\be\label{e2.12}
\calE :=\|u_{a,\delta}(t)-y\| \leq \|u_{a,\delta}(t)-u_a(t)\|+\|u_a(t)-y\|. \ee
One has
\be\label{e2.13}
\|u_{a,\delta}(t)-u_a(t)\| \leq \frac{\|f_\delta-f\|}{a}\| I-e^{iAt-at}\| \leq \frac{2\delta}{a}, \ee
and
\be\label{e2.14}
\|u_a(t)-y\| \leq \frac{e^{-at}}{a}\| f \|+\eta(a),\ee
where $\eta(a)$ is defined in (2.8), $\lim_{a\to 0}\eta(a)=0$.
Since $\|f\| \leq \|f_\delta\|+\delta\leq c$,  one obtains from \eqref{e2.12}-\eqref{e2.14}:
\be\label{e2.15}  \lim_{\delta\to 0}\calE=0 \ee
provided that $t=t_\delta$, $a=a(\delta)$ and 
\be\label{e2.16}
\lim_{\delta\to 0}t_\delta=\infty, \qquad \lim_{\delta\to 0} 
a(\delta)=0, 
\qquad \lim_{\delta\to 0}\frac{e^{-a(\delta) t_\delta}}{a(\delta)}=0,
\qquad \lim_{\delta\to 0} \frac{\delta}{a(\delta)}=0. \ee

\refT{2} is proved.
\end{proof}

\begin{remark}\label{R:2}
There are many choices of $t_\delta$ and $a(\delta)$ satisfying relations \eqref{e2.16}.
If one has an estimate of the rate of decay of $\eta(a)$ as $a\to 0$, then one may obtain
some rate of convergence of $\calE$ to zero as $\delta\to 0$.
However, it is impossible, in general, to get a rate of decay of $\eta(a)$ as $a\to 0$
without additional assumptions on the data $f$ or on the solution $y$. A typical assumption
is $y=Az$, that is, $y\in R(A)$. If fractional powers of $A$ are defined (which is the case 
when  $A\geq 0$, for example) then one may assume $y=A^\gamma z$, 
$\gamma>0$. 
Let us show how to get the rate of decay of $\eta(a)$ under such assumptions. 
Assume, for example, that $y=Az$. Then
\be\label{e2.17}
\eta^2(a)=\int^\infty_{-\infty} \frac{a^2 s^2}{a^2+s^2} d(E_sz,z)\leq a^2\|z\|^2,\ee
and the error bound is
\be\label{e2.18}
\calE\leq\frac{2\delta}{a} + c_1\frac{e^{-at}}{a}+c_2 a,
\qquad c_2=\|z\|, \quad c_1=\|f_\delta\|+\delta. \ee
Choose, for example, $a=\delta^\gamma$, $0<\gamma<1$, 
and $t_\delta=\delta^{-\mu}$, $\mu>\gamma$. 
Then $\lim_{\delta\to 0} \frac{e^{-a(\delta)t_\delta}}{a(\delta)}=0$ 
and \eqref{e2.15} holds with the rate $\delta^\nu$, $\nu=min(1-\gamma,\gamma)$.
If $\gamma=\frac{1}{2}$ then $\max_{0<\gamma<1}\min(1-\gamma,\gamma)$ is 
equal to $\frac{1}{2}$,
and for $\gamma=\frac{1}{2}$ one gets $\calE=O(\delta^{1/2})$ if
$t_\delta=\delta^{-\mu}$, $\mu>\frac{1}{2}$.
\end{remark}

\section{Second version of the DSM}\label{S:3}
Consider problem \eqref{e2.1} with $a=a(t)$.
Let us assume that $a(t)>0$ is a continuous function monotonically decaying to zero as
$t\to\infty$ and
\be\label{e3.1}
0<a(t)\searrow 0,  \qquad a'+a^2\in L^1(0,\infty),
\qquad \int^\infty_0 a(s)ds=\infty. \ee
The solution to this problem is
\be\label{e3.2}
u(t)=\int^t_0 e^{iA(t-s)-\int^t_s a(p)dp} ds(-if).\ee

\begin{theorem}\label{T:3}
Under the above assumptions one has
\be\label{e3.3}   \lim_{t\to\infty}\|u(t)-y\|=0. \ee
\end{theorem}

\begin{proof} [Proof of Theorem 3.]
Since $f=Ay$, one  gets, integrating by parts, 
\bee  u(t)=e^{iAt} e^{-iAs-\int^t_s adp}y|_0^t-
\int^t_0 e^{iA(t-s)} a(s) e^{-\int^t_{s} a(p)dp} ds y.\eee
Thus
\be\label{e3.4}
u(t)=y-e^{iAt-\int^t_0 adp} y -\int^t_0 e^{iA(t-s)} a(s) e^{-\int^t_s adp} 
dsy, \ee
and
\be\label{e3.5}
\|u(t)-y\|\leq e^{-\int^t_0 adp} \|y\| + 
\|\int^t_0 e^{iA(t-s)}a(s) e^{-\int^t_s adp}ds y\|:=J_1+J_2.\ee
By the last assumption \eqref{e3.1} one gets
\be\label{e3.6}  \lim_{t\to\infty}J_1=0  \ee

We now prove that
\be\label{e3.7}  \lim_{t\to\infty}J_2=0. \ee
Using the spectral theorem, one gets
\be\label{e3.8}
J^2_2 =\int^\infty_{-\infty} d(E_\lambda y,y)
\bigg| \int^t_0 e^{i\lambda(t-s)} a(s)e^{-\int^t_s adp} ds \bigg|^2.  \ee
Let us prove that
\be\label{e3.9}
\lim_{t\to\infty} \int^t_0 e^{i\lambda(t-s)} a(s) e^{-\int^t_s adp}ds=0,
\qquad \forall \lambda\not=0.\ee
From \eqref{e3.8}, \eqref{e3.9} and the assumption $y\perp N$,
the conclusion \eqref{e3.7} follows.

Let us verify \eqref{e3.9}. Integrating by parts one gets:
\be\label{e3.10}
J_3:=\int^t_0 e^{i\lambda(t-s)} a(s) e^{-\int^t_s adp}ds
=\frac{e^{i\lambda(t-s)}}{-i\lambda} a(s) e^{-\int^t_s adp} 
\bigg|^t_0 +\frac{1}{i\lambda} \int^t_0 
e^{i\lambda(t-s)} [a'(s)+a^2(s)]e^{-\int^t_sa(p)dp}ds. \ee
Thus
\be\label{e3.11}
J_3= \frac{a(t)}{-i\lambda} + \frac{e^{i\lambda t-\int^t_0 adp}a(0)}{i\lambda} + J_4,
\qquad  \lambda\not= 0,\ee
where $J_4$ denotes the last integral in \eqref{e3.10}.
The first two terms in \eqref{e3.11} tend to zero as $t\to\infty$ because of the 
assumptions about $a(t)$.

Assumptions \eqref{e3.1} imply that
\be\label{e3.12}  \lim_{t\to\infty} J_4=0. \ee
Thus, \refT{3} is proved.
\end{proof}

\begin{remark} \label{R:3}
For example, the function $a(t)=\frac{c_0}{(c_1+t)^b}$, where $c_0,c_1>0$
are arbitrary constants and $b\in(\frac{1}{2},1)$ is a constant, satisfies
assumptions \eqref{e3.1}.

It is interesting to check numerically the efficiency of the algorithm based on \refT{3}.
\end{remark}

Let us prove that \refT{3} yields a stable solution to equation \eqref{e1.1}.

\begin{theorem}\label{T:4}
There exists a stopping time $t_\delta$, $\lim_{\delta\to 0} 
t_\delta=\infty$,
such that \eqref{e1.6} holds with $u_\delta=u_\delta(t_\delta)$,
where $u_\delta(t)$ is the solution to problem \eqref{e2.1} with 
$a=a(t)$ and $f_\delta$ in place of $f$, $\|f_\delta-f\|\leq\delta$.
\end{theorem}

\begin{proof}[Proof of \refT{4}.]
One has
\be\label{e3.13}
\|u_\delta(t)-y\| \leq \|u_\delta(t)-u(t)\| + \|u(t)-y\|, \ee
where $u(t)$ solves problem \eqref{e2.1} with $a=a(t)$ and exact data. 
We have proved in \refT{3} that
\be\label{e3.14}  \lim_{t\to\infty} \|u(t)-y\|=0.  \ee
We have
\be\label{e3.15}
\|u_\delta(t)-u(t)\| \leq \int^t_0 e^{-\int^t_s a(p)dp}ds
\|f_\delta -f\| \leq \frac{\delta}{a(t)}.\ee
Here the estimate $\int^t_0 e^{-\int^t_s adp}ds \leq\frac{1}{a(t)}$ was used.
This estimate is derived easily:
\bee
\int^t_0 e^{-\int^t_s adp}ds
\leq\frac{1}{a(t)}\int^t_0 a(s) e^{-\int^t_s adp}dp
=\frac{1}{a(t)} e^{-\int^t_s adp}\bigg|^t_0
=\frac{1}{a(t)}\left(1-e^{-\int^t_0 adp}\right) \leq \frac 1 {a(t)}. \eee
Choose $t_\delta$ so that
\be\label{e3.16}
\lim_{\delta\to 0} t_\delta=\infty, \qquad \lim_{\delta\to 0} \frac{\delta}{a(t_\delta)}=0. \ee
This is obviously possible to do.

Then \eqref{e3.13}-\eqref{e3.15} imply
\be\label{e3.17}
\lim_{\delta\to 0} \|u_\delta(t_\delta)-y\|=0. \ee

\refT{4} is proved.
\end{proof}

\section{Third version of DSM for solving equation (1.1)} \label{S:4}

\be\label{e4.1}
\dotu =-u+T^{-1}_{a(t)} A^\ast f, \qquad u(0)=0, \ee
where $f=Ay$, $y\perp N$, and $a(t)>0$ is a monotonically
decaying continuous function such that $\lim_{t\to 
\infty}a(t)=0$ and
$\int^\infty_0 a(t)dt=\infty$.
We could take the initial condition $u(0)=u_0$ arbitrary. The contribution 
to the solution to problem \eqref{e4.1} which comes from the initial
condition $u_0$ is the term $u_0e^{-t}$. It decays exponentially fast and 
our arguments do not depend on this term essentially. To simplify and 
shorten our argument we take $u_0=0$.

\begin{theorem}\label{T:5}
$\ds \lim_{t\to\infty} u(t)=y$.
\end{theorem}

\begin{proof}[Proof of Theorem 5.]
One has
\bee u(t)=\int^t_0 e^{-(t-s)}  T^{-1}_{a(s)} A^\ast Ayds
=\int^t_0 e^{-(t-s)} yds-\int^t_0 e^{-(t-s)} a(s) T^{-1}_{a(s)} yds. \eee
Thus
\be\label{e4.2}
\|u(t)-y\| \leq e^{-t} \|y\| + \int^t_0 e^{-(t-s)}
a(s)\|T^{-1}_{a(s)} y\| ds. \ee
One can easily check that if $b(s)$  is a continuous function on $[0,\infty)$ and 
$b(\infty)=\lim_{s\to\infty}b(s)$ exists, then
\be\label{e4.3}
\lim_{t\to\infty} \int^t_0 e^{-(t-s)} b(s)ds=b(\infty). \ee
Thus, \refT{5} will be proved if one checks that
\be\label{e4.4}
\lim_{a\to 0} a\|T^{-1}_a y\|=0 \qquad \forall  y\in H.\ee
To prove \eqref{e4.4} one writes, using the spectral theorem:
\be\label{e4.5}
a^2 \|T^{-1}_ay\|^2= \int^\infty_0 \frac{a^2}{(s+a)^2} d(E_s y,y)
=\|P_N y\|^2=0.\ee
\refT{5} is proved.
\end{proof}

Let us prove that the DSM method \eqref{e4.1} yields a stable solution
to problem \eqref{e1.1}. 
Assume that $f$ is replaced by $f_\delta$, $\|f_\delta-f\|\leq\delta$,
in equation \eqref{e4.1} and denote by $u_\delta(t)$ the corresponding 
solution. Then, using the estimate \eqref{e3.13}, one gets
\be\label{e4.6}
\lim_{\delta\to 0} \|u_\delta(t_\delta)-y\|=0, \ee
provided that
\be\label{e4.7}
\lim_{\delta\to 0} t_\delta=\infty, \qquad
\lim_{\delta\to 0} \frac{\delta}{\sqrt{a(t_\delta)}}=0. \ee
To check the sufficiency of the second condition \eqref{e4.7} for \eqref{e4.6} 
to hold, one proceeds as follows:
\be\label{e4.8}
\begin{aligned}
\|u_\delta(t)-u(t)\| 
   &= \| \int^t_0 e^{-(t-s)} T^{-1}_{a(s)} A^\ast (f_\delta-f)ds\| \\
   &\quad \leq \delta \int^t_0 e^{-(t-s)} \frac{1}{2\sqrt{a(s)}} ds
   \leq \frac{\delta}{2\sqrt{a(t)}}.
   \end{aligned} \ee
Here we have used the monotonicity of $a(t)$, which implies $a(t)\leq a(s)$ 
if $t\geq s$, and the estimate $\|T^{-1}_a A^\ast\|\leq\frac{1}{2\sqrt{a}}$,
which was proved earlier.

Let us state the result we have proved.

\begin{theorem}\label{T:6}
If $t_\delta$ is chosen so that \eqref{e4.7}  holds, then the solution $u_\delta(t)$ to 
problem \eqref{e4.1} with noisy data $f_\delta$ in place of $f$ satisfies \eqref{e4.6}.
\end{theorem}

\section{A new discrepancy principle}\label{S:5}

The usual discrepancy principle is described in \refS{1}. It requires solving nonlinear equation
$\|Au_{a,\delta}-f\|=C\delta$, $C=const$, $1<C<2$,
where $u_{a,\delta}=T^{-1}_a A^\ast f_\delta$.
Thus one has to know the exact minimizer $u_{a,\delta}$ of the functional 
$F(u)=\|Au-f_\delta\|^2 +a\|u\|^2$,
or the exact solution of the equation $T_a u=A^\ast f_\delta$.

In this Section we discuss the following question:

\textit{How does one formulate the discrepancy principle in the case when
$u_{a,\delta}$ is not the exact solution of the minimization problem
$F(u)=\min$, but an approximate solution.}

Let us state the result.

\begin{theorem}\label{T:7}
Assume that $A$ is a bounded linear operator in a Hilbert space $H$,
that $f=Ay$, $y\perp N$, $\|f_\delta-f\|\leq\delta$, $\|f_\delta\|>C\delta$,
$C=const$, $C\in(1,2)$, and $u_{a,\delta}$ is any element 
which satisfies the inequality:
\be\label{e5.1}
F(u_{a,\delta})\leq m+(C^2-1-b)\delta^2, \ee
where
\be\label{e5.2}
F(u)=\|Au-f_\delta\|^2 + a\|u\|^2, \qquad m=\inf_{u\in H} F(u), \ee
$b=const>0$ and $C^2>1+b$.

Then equation
\be\label{e5.3}  \|Au_{a,\delta}-f_\delta\|=C\delta \ee
has a solution for any fixed $\delta>0$, $\lim_{\delta\to 0} a(\delta)=0$, and 
\be\label{e5.4}  \lim_{\delta\to 0} \|u_\delta-y\|=0 \ee
where $u_\delta=u_{a(\delta),\delta}$ and $a(\delta)$ solves \eqref{e5.3}.
\end{theorem}

\begin{proof}[Proof of \refT{7}]
To prove the existence of a solution to \eqref{e5.3} denote 
$\|Au_{a,\delta}-f_\delta\|$  by $h(\delta,a)$, check that $h(\delta,+0)<C\delta$,
$h(\delta,\infty)>C\delta$ and note that $h(\delta,a)$ is a continuous function of $a$
on the interval $(0,\infty)$. Then equation \eqref{e5.3} has a solution $a=a(\delta)$.
As $a\to\infty$ one has
\bee a\|u_{a,\delta}\|^2  \leq F(u_{a,\delta})\leq m+(c^2-1-b)\delta^2
\leq F(0)+(C^2-1-b)\delta^2 \eee
so, with $c:=F(0)+(C^2-1-b)\delta^2$, one obtains
\bee \|u_{a,\delta}\|\leq \frac{c}{\sqrt{a}}, \qquad a\to\infty. \eee
Thus 
\be\label{e5.5}
h(\delta,\infty)=\|A0-f_\delta\| = \|f_\delta\|>C\delta. \ee
As $a\to 0,$ one has
\bee h^2(\delta,a)\leq F(u_{a,\delta}) 
\leq m+(C^2-1-b)\delta^2 \leq F(y)+(C^2-1-b)\delta^2. \eee
One has $F(y)=\delta^2+a\|y\|^2$, so
$h^2(\delta,a)\leq (C^2-b)\delta^2+a\|y\|^2$, and
\be\label{e5.6}
h(\delta,+0)\leq (C^2-b)^{\frac{1}{2}}\delta < C\delta. \ee
Finally, the continuity of $h(\delta,a)$ with respect to  $a$ for any fixed $\delta>0$
follows from the continuity of the bounded operator $a$ 
and the assumed continuity of $u_{a,\delta}$ with respect to $a$.
Thus, the existence of a solution $a=a(\delta)>0$ of equation \eqref{e5.3}
is proved.
One takes a solution for which  $\lim_{\delta\to 0}a(\delta)$.
Such solution exists because $h(\delta,+0)$ and $m=m(\delta,a)$ tend to zero
as $\delta\to 0$ and $a\to 0$.

Let us prove \eqref{e5.4}. One has 
\be\label{e5.7}
F(u_\delta)=\|Au_\delta-f_\delta\|^2
+ a(\delta)\|u_\delta\|^2
\leq \|Ay-f_\delta\|^2
+a(\delta) \|y\|^2
\leq \delta^2+a(\delta)\|y\|^2. \ee
Since $\|Au_\delta-f_\delta\|=C\delta$, $C>1$ it follows from \eqref{e5.7}
that
\be\label{e5.8}  \|u_\delta\|\leq\|y\|. \ee
Thus one may assume that $u_\delta\rightharpoonup u$ as $\delta\to 0$.

Let us prove that $Au=f$. First we observe that
\bee \|Au_\delta-f\| \leq \|Au_\delta-f_\delta\| + \|f_\delta-f\| \leq C\delta+\delta, \eee
so
\be\label{e5.9}   \lim_{\delta\to 0} \|Au_\delta-f\|=0  \ee
Secondly, for any $v\in H$ we have
\be\label{e5.10}
(f-Au,v)=\lim_{\delta\to 0} (Au_\delta-Au,v)
=\lim_{\delta\to 0} (u_\delta-u, A^\ast v)=0, \ee
because $u_\delta\rightharpoonup u$. Since $v$ is arbitrary, one concludes
from \eqref{e5.10} that $Au=f$. From \eqref{e5.8} it follows that 
$\|u\|\leq \|y\|$.  Since the minimal-norm solution to equation $Au=f$ is
unique, one obtains $u=y$. Thus $u_\delta\rightharpoonup y$, $\|u_\delta\|\leq\|y\|$.
This implies \eqref{e5.4} as follows from \eqref{e1.18}.

\refT{7} is proved.
\end{proof}

\section{Discrepancy principle does not yield convergence uniformly
with respect to the data}\label{S:6}

In this Section we make the following assumption.

{\bf Assumption} A): {\it  $A$ is a linear bounded operator in a Hilbert  
space 
$H$, $N:=N(A)=N(A^\ast):=N^\ast=\{0\}$, $A^{-1}$ is unbounded, $Ay=f$, 
 $\|f_\delta-f\|\leq\delta$, $\|f_\delta\|>\delta$}.
 
 Let $a=a(\delta)$ be chosen by the discrepancy principle,
 \be\label{e6.1}
 \|Au_{a,\delta}-f_\delta\|=C\delta,  \qquad u_\delta=u_{a(\delta),\delta} 
 = T^{-1}_{a(\delta)} A^\ast f_\delta. \ee
 Consider the set $S_\delta$:
 \bee S_\delta:=\{v:\|Av-f_\delta\|\leq\delta\}. \eee
We are interested in the following question: 
{\it given $\{f_\delta,\delta\}_{\delta\in(0,\delta_0)}$, 
where $\delta_0>0$ is a small number, and  assuming that $a(\delta)$ is
the solution to \eqref{e6.1}, can one guarantee uniformity with respect to the data $f$
convergence?} 

In other words, is it true that
\be\label{e6.2}
\lim_{\delta\to 0} \sup_{v\in S_\delta} \|u_\delta-v\|=0. \ee
The answer is no.

\begin{theorem} \label{T:8}
There exist $f_\delta$ such that
\be\label{e6.3}
\lim_{\delta\to 0} \sup_{v\in S_\delta} \|u_\delta-v\|\geq c>0,
\qquad c=const. \ee
\end{theorem}

\begin{proof} [Proof of Theorem 8.]
Denote $T^{-1}_a A^\ast:=G$, $u_\delta=G f_\delta$, 
$\|G\|=\frac{1}{2\sqrt{a}}$. 
We have proved in \eqref{e1.10} that $\|G\|\leq\frac{1}{2\sqrt{a}}$,
but the equality sign holds because in \eqref{e1.10} $U$ is unitary
under {\bf Assumptions} A). Thus, one can find an element $p=p_a$,
$\|p\|=\frac{\delta}{2}$, such that
\be\label{e6.4}
\|Gp\|\geq \frac{1}{2} \|G\|\ \|p\|=\frac{\delta}{8\sqrt{a}}. \ee
Assumptions A) imply that the ranges $R(A)$ and $R(T)$ are dense in $H$. 
Thus one can find an element $z=z_{a,\delta}$ such that
\be\label{e6.5}
\|f_\delta-AT^bz-p\|\leq \frac{\delta}{8}, \qquad b\in(0,1), \ b=const. \ee
For any $v$ one has
\be\label{e6.6}
\|Gf_\delta-v\| \leq \|Gf_\delta-G Av\| + \|GAv-v\|. \ee
Take $v=T^bz$ and let $M>0$ be an arbitrary large fixed $a_n$ constant. 
Then
\be\label{e6.7}
\lim_{\delta\to 0} \sup_{v\in S_\delta, v=B^bz, \|z\|\leq m}
\|GAv-v\|=0, \ee
because
\be\label{e6.8}
\|GAv-v\| = \|T^{-1}_a Tv-v\| = \|-aT^{-1}T^bz\|
=a \sup_{s\geq 0} \frac{s^b}{s+a}=ca^b, \ee
where $c=b^b(1-b)^{1-b}$.

From \eqref{e6.7} and \eqref{e6.6} one sees that
\be\label{e6.9}
\lim_{\delta\to 0} \sup_{v\in S_\delta, v=T^bz, \|z\|\leq M}
\|Gf_\delta-v\|=0.\ee
if and only if
\be\label{e6.10}
\lim_{\delta\to 0} \sup_{v\in S_\delta, v=T^bz, \|z\|\leq M} 
\|Gf_\delta-GAv\|=0. \ee

Take $z=z_a$, $v=T^bz$. Then
\be\label{e6.11}
\|f_\delta -Av\| \leq \|p\|+\|f_\delta-AB^bz-p\|
\leq \frac{\delta}{2}+\frac{\delta}{8}=\frac{5\delta}{8}\leq \delta, \ee
and, using \eqref{e6.5}, one gets:
\be\label{e6.12}
\|Gf_\delta-GAv\| = \|G(f_\delta-Av-p)+Gp \| \geq \|Gp\| - 
\|G\|\frac{\delta}{8}
\geq \frac{\delta}{\sqrt{a}} 
\left(\frac{1}{8}-\frac{1}{16}\right)=\frac{\delta}{16\sqrt{a}}. \ee
If $\frac{\delta}{\sqrt{a}}\geq c>0$, then, according 
to \eqref{e6.12},  \eqref{e6.10} fails. 

Let us find $f_\delta$ such that for $a=a(\delta)$, defined by the 
discrepancy
principle, one has $\frac{\delta}{\sqrt{a}}\geq c>0$. This will 
complete the proof of \refT{8}. Let us assume for simplicity that 
$A=A^\ast>0$ is compact. Then $T=A^\ast A=A^2$, and
equation \eqref{e6.1} becomes
\be\label{e6.13}
C^2\delta^2
  = \|[A(A^2+a)^{-1}A-I]f_\delta\|^2
  = \sum^\infty_{j=1} \left[\frac{\lambda^2_j}{\lambda^2_j+a} -1\right]^2 
  |f_{\delta j}|^2
  =\sum^\infty_{j=1} \frac{a^2|f_{\delta j}|^2}{(\lambda_j+a)^2}. \ee
Here $\lambda_j$ are the eigenvalues of $A^2$, $f_{\delta 
j}=(f_\delta,\varphi_j)$,
$A^2\varphi_j=\lambda_j\varphi_j$, $\|\varphi_j\|=1$.
Assume, for example, that $\lambda_j=\frac{1}{j}$ and 
$|f_{\delta j}|^2=\frac{1}{j^2}$. Then \eqref{e6.13} becomes
\be\label{e6.14}
C^2\frac{\delta^2}{a^2} = \sum^\infty_{j=1} \frac{j^{-2}}{(j^{-1}+a)^2}. \ee
Note that
\bee\begin{aligned}
I_1:&=\int^\infty_1 \frac{x^{-2}}{(x^{-1}+a)^2}dx\\
&=\int^1_0 \frac{ds}{(s+a)^2}
=-(s+a)^{-1}\bigg|^1_0 
=a^{-1}[1-\frac{a}{a+2}]
=\frac{1}{a}[1+O(a)],
\quad a\to 0\end{aligned}\eee

This and \eqref{e6.14} imply $\frac{\delta}{\sqrt{a}}\geq c>0$ 
as $\delta\to 0$.
\refT{8} is proved.
\end{proof}

\section{Iterative processes for solving equation \eqref{e1.1}}\label{S:7}

In this Section convergent iterative processes for solving equation \eqref{e1.1} 
are constructed in the case when $A$ is a closed, densely defined in $H$, 
unbounded operator.
Consider the process
\be\label{e7.1} 
u_{n+1}=Bu_n+T^{-1}_a A^\ast f, 
\quad u_1=u_1, \quad u_1\perp N, \quad B:=aT^{-1}_a,\ee
where $a=const>0$ and the initial element $u_1$ is arbitrary in the
subspace $N^\perp$, $N:=N(A)=N(T)$, $T=A^\ast A$, $T_a=T+aI$. Note that 
$B\geq 0$, $\|B\|\leq 1$.

\begin{theorem}\label{T:9}
Under the above assumptions one has
\be\label{e7.2}   \lim_{n\to\infty} \|u_n-y\|=0. \ee
\end{theorem}

\begin{proof}[Proof of \refT{9}]
Let $w_n=u_n-y$. Then 
\be\label{37.3}
w_{n+1}=Bw_n=B^n w, \qquad w:=u_1-y, \qquad w\perp N. \ee
Let us prove that
\be\label{e7.4}   \lim_{n\to\infty} \|B^nw\|=0. \ee
If \eqref{e7.4} is verified, then \refT{9} is proved. We have
\be\label{e7.5}
\|B^nw\|^2 = \int^\infty_0 \frac{a^{2n}}{(a+s)^{2n}}d(E_sw,w)
=\int_{s>b}+ \int_{0\leq s\leq b}:= J_1+J_2, \ee
where $E_s$ is the resolution of the identity corresponding to the operator
$T\geq 0$, and $b>0$ is a small number which will be chosen later.
For any fixed $b>0$ one has $\lim_{n\to\infty}J_1=0$
because $\frac{a}{a+s}\leq\frac{a}{a+b}<1$ if $s\geq b$.
On the other hand, $J_2\leq \int^b_0d(E_s w,w)$,
and $\lim_{b\to 0} \int^b_0 d(E_s w,w)=0$ because $w\perp N=E_0H$.
Therefore, given an arbitrary small number $\eta>0$ one can choose
$b>0$ such that $J_2\leq\eta/2$. Fix this $b$ and choose $n$ sufficiently
large so that $J_1\leq \eta/2$. Then $\|B^nw\|^2\leq\eta$.
Since $\eta$ is arbitrarily small, we have proved \eqref{e7.4}
\refT{9} is proved.
\end{proof}

\begin{remark}\label{R:4}
Iterative process \eqref{e7.1} yields a stable solution of equation\eqref{e1.1}.
Indeed, let $f_\delta$ be given, $\|f_\delta-f\|\leq\delta$, and
let $u_{n,\delta}$ be defined by \eqref{e7.1} with $f_\delta$ in place of $f$.
Let $w_{n,\delta}=u_{n,\delta}-y$. Then
\be\label{e7.6}
\begin{aligned}
w_{n_1,\delta} 
  & =Bw_{n,\delta} + T^{-1}_a A^\ast(f_\delta-f), \hbox{\ \ so\ \ }\\
 w_{{n+1},\delta}
  &= \sum^n_{j=1} B^j T^{-1} A^\ast(f_\delta-f)+B^n(u_1-y), 
 \qquad (u_1-y)\perp N.\end{aligned}\ee

We have proved above that
\be\label{e7.7} \|B^n(u_1-y)\| :=\calE(n)\to 0 \hbox { as } n\to \infty. 
\ee
One has:
\be\label{e7.8}
 \bigg \|\sum^n_{j=0} B^j T^{-1}_a A^\ast (f_\delta-f)\bigg\|
  \leq\frac{(n+1)\delta}{2\sqrt{a}},\ee
because $\|B\|\leq 1$ and $\|T^{-1}_a A^\ast\|\leq\frac{1}{2\sqrt{a}}$.
From \eqref{e7.7} and \eqref{e7.8} one finds the stopping rule, i.e.,
the number $n(\delta)$ such that 
$\lim_{\delta\to 0} \|w_{n(\delta),\delta} \|=0$.
This $n(\delta)$ is found for any fixed small $\delta$ as the
minimizer for the problem
\be\label{e7.9}
\frac{(n+1)\delta}{2\sqrt{a}} + \calE(n)=min. \ee
Alternatively, one can find $n_1(\delta)$ from the equation
\be\label{e7.10}
\calE(n)=\frac{(n+1)\delta}{2\sqrt{a}}.\ee
Clearly $n(\delta)$ and $n_1(\delta)$ tend to $\infty$ as 
$\delta\to 0$.
\end{remark}

\section{Discrepancy principle for DSM}\label{S:8}
In this Section we formulate and justify a discrepancy principle for DSM.

Let us start with the version \eqref{e4.1}. We assume that $a(t)>0$
is a monotonically decaying twice continuously differentiable 
function, $\lim_{t \to \infty}[a(t)+|\dot{a}|+\ddot{a}]=0$, $\ddot{a}>0$,
$\lim_{t \to \infty}\frac {\dot a(t)}{a(t)}=0$, for example, $a(t)=\frac 
{c_1}{(c_0+t)^b}$, where $c_1,\,c_0$ and $b$ are 
positive constants, $b\in (0.5, 1)$. For this $a(t)$ all the assumptions
(3.1) hold.

\begin{theorem}\label{T:10}

The equation
\be\label{e8.1}
\|A T^{-1}_{a(t)} A^\ast f_\delta-f_\delta\| =C\delta
\qquad C=const, \quad 1<C<2, \ee
has a solution $t=t_\delta$, $\lim_{\delta\to 0} t_\delta=\infty$, 
such that \eqref{e4.6} holds, where $u_\delta(t)$ is the solution to
\eqref{e4.1} with $f_\delta$ in place of $f$, $\|f_\delta\|>C\delta.$
\end{theorem}

\begin{proof}[Proof of \refT{10}]
We have proved earlier that equation \eqref{e8.1} has a unique solution
$a=a_\delta$ and  $\lim_{\delta\to 0} a_\delta=0$.
If $a(t)$ is a monotonically decaying function such that 
$\lim_{t\to\infty}a(t)=0$, then the equation $a_\delta=a(t)$
defines uniquely $t=t_\delta$, such that $a(t_\delta)=a_\delta$,
and $\lim_{\delta\to 0} t_\delta=\infty$.

Let us sketch the proof of relation \eqref{e4.6}, where 
 $u_\delta(t_\delta)=\int^{t_\delta}_0 e^{-(t_\delta-s)} T^{-1}_{a(s)} 
A^\ast f_\delta ds$ and $t_\delta\to\infty$ as $\delta\to 0$.
We have proved (cf (1.17)) that $||T^{-1}_{a(t_\delta)} A^\ast 
f_\delta||\leq ||y||$, and 
\be\label{e8.2}
\lim_{\delta\to 0}||T^{-1}_{a(t_\delta)} A^\ast f_\delta -y||=0.
\ee 
(cf (1.18)). It is clear that
$\lim_{t \to \infty}\int_0^t e^{-(t-s)} g(s)ds=g(\infty)$ provided that 
$g$ is a continuous function and there exists $g(\infty):=\lim_{t \to 
\infty}g(t)$. 

Note that $\lim_{s\to t_\delta}||T^{-1}_{a(t_\delta)}A^\ast f_\delta-
T^{-1}_{a(s)}A^\ast f_\delta||=0$. We have 
\be\label{e8.3}
\int_0^{t_\delta}e^{-(t_\delta -s)}T^{-1}_{a(s)}A^\ast f_\delta 
ds=y+o(1) \hbox { as } \delta \to 0.\ee
To check this we use equation (8.2), the formula
$T^{-1}_{a(t_\delta)}-T^{-1}_{a(s)}=T^{-1}_{a(s)}[a(t_\delta)-a(s)]
T^{-1}_{a(t_\delta)}$, the estimates $||T^{-1}_{a(t_\delta)}A^\ast 
f_\delta||\leq ||y||$ and $||T^{-1}_{a(s)}||\leq \frac 1 {a(s)}$, and the 
relation
$\lim_{t \to \infty}\int_0^t e^{-(t-s)} \frac {a(s)-a(t)}{a(s)}ds=0$,
which holds due to our assumptions on $a(t)$. Theorem 10 is 
proved.
\end{proof}


\end{document}